\documentclass{article}
\usepackage{amsmath}
\usepackage{amssymb}
\usepackage{amsthm}
\newtheorem{theorem}{Theorem}
\newtheorem{lemma}[theorem]{Lemma}
\newtheorem{corollary}[theorem]{Corollary}
\newtheorem{definition}[theorem]{Definition}

\newtheorem{sublemma}[theorem]{Sublemma}
\newtheorem{proposition}[theorem]{Proposition}
\newtheorem{desired result}[theorem]{Desired result}
\newtheorem{conjecture}[theorem]{Conjecture}
\newtheorem{question}[theorem]{Question}
\newtheorem{notation}[theorem]{Notation}

\newcommand{\begd}{\begin{displaystyle}}
\newcommand{\gl}{\lambda}

\newcommand{\ga}{\alpha}
\newcommand{\gb}{\beta}
\newcommand{\gd}{\delta}
\newcommand{\gD}{\Delta}
\newcommand{\gs}{\sigma}
\newcommand{\tr}{\textrm{tr}\,}
\newcommand{\ve}{\varepsilon}
\newcommand{\om}{\omega}
\newcommand{\mbq}{\mathbb{Q}}

\newcommand{\mbz}{\mathbb{Z}}

\newcommand{\gal}{\textrm{Gal}\,}
\newcommand{\ol}[1]{\overline{#1}}
\newcommand{\mc}[1]{\mathcal{#1}}

\title{Averages of elliptic curve constants}
\author{Nathan Jones}
\date{jones@dms.umontreal.ca}

\begin{document}

\maketitle
\begin{abstract}
We compute the averages over elliptic curves of the constants occurring in the Lang-Trotter conjecture, the Koblitz conjecture, and the cyclicity conjecture.  The results obtained confirm the consistency of these conjectures with the corresponding ``theorems on average'' obtained recently by various authors.
\end{abstract}

\section{Introduction} \label{introduction}

Let $E$ be an elliptic curve defined over the rational numbers and for a prime $p$ of good reduction for $E$, let $E_p$ denote the reduction of $E$ modulo $p$.  There are various conjectured asymptotics for functions which count good primes $p$ up to $x$ for which the reduced curve $E_p$ has certain properties.  In this paper we will focus on three  such questions, although our methods are applicable to a wider range of problems.  For a fixed integer $r$, let
\[
\pi_{E,r}(x) = |\{ p \leq x : p \nmid \gD_E, \; a_p(E) = r \}|,
\]
where $a_p(E) = p+1 - |E_p(\mbz/p\mbz)|$ is the trace of the Frobenius endomorphism of $E$ at $p$.  Lang and Trotter \cite{lt}, using a probabilistic model consistent with the Chebotarev density theorem and the Sato-Tate conjecture, predicted an asymptotic for $\pi_{E,r}(x)$:
\begin{conjecture}  \label{langtrotter} (Lang-Trotter)  Assume that either $E$ has no complex multiplication or that $r \neq 0$.  Then 
\[
\pi_{E,r}(x) \sim C_{E,r} \frac{\sqrt{x}}{\log x} \quad\quad \textrm{ as } \; x \rightarrow \infty,
\]
where $C_{E,r}$ is a specific constant.  
\end{conjecture}
We will describe the constant $C_{E,r}$ in detail in section \ref{theconstants}.  The second conjecture we will consider involves the counting function
\[
\pi_{E, \text{prime}}(x) := | \{ p \leq x : p \nmid \gD_E, \; |E(\mbz/p\mbz)| \text{ is prime } \} |.
\]
\begin{conjecture} \label{koblitz} (Koblitz)  
\[
\pi_{E, \text{prime}}(x) \sim C_{E, \text{prime}} \frac{x}{(\log x)^2} \quad\quad \textrm{ as } \; x \rightarrow \infty,
\]
where $C_{E, \text{prime}}$ is a specific constant.
\end{conjecture}
Finally we will consider the cyclicity conjecture, which has been settled conditionally by Serre \cite{serrecyclic} and unconditionally in the CM case by Murty \cite{murty} and also by Cojocaru \cite{cojocaru}.  Let
\[
\pi_{E, \text{cyclic}}(x) := | \{ p \leq x : p \nmid \gD_E, \; E(\mbz/p\mbz) \text{ is a cyclic group } \} |.
\]
\begin{conjecture} \label{cyclicity} (Cyclicity conjecture)
\[
\pi_{E, \text{cyclic}}(x) \sim C_{E, \text{cyclic}} \frac{x}{\log x} \quad\quad \textrm{ as } \; x \rightarrow \infty,
\]
where $C_{E, \text{cyclic}}$ is a specific constant.
\end{conjecture}
Recently, various authors have proven that these conjectures ``hold on average over elliptic curves.''
More precisely, for parameters $A = A(x)$ and $B = B(x)$, let $\mc{C} = \mc{C}(x)$ denote the set of elliptic curves $Y^2 = X^3 + aX + b$ with $(a,b) \in \left( [-A,A]\times[-B,B] \right) \cap \mbz^2$.
Fouvry and Murty \cite{fm} (in case $r = 0$) and later David and Pappalardi \cite{dp} (in case $r \neq 0$) proved Conjecture \ref{langtrotter} on average:  for any $\ve > 0$, if $\min \{ A(x), B(x) \} \geq x^{1+\ve}$ then
\begin{equation} \label{langtrotterthm}
\frac{1}{|\mc{C}(x)|} \sum_{{\begin{substack} {E \in \mc{C}(x)} \end{substack}}} \pi_{E, r}(x) \sim C_r \frac{\sqrt{x}}{\log x}, \quad \textrm{as } x \rightarrow \infty,
\end{equation}
where $C_r$ is a specific constant.  S. Baier \cite{baier} has recently shortened the length of the average, replacing ``$x^{1+\ve}$'' with ``$x^{3/4+\ve}$.''

Balog, Cojocaru and David \cite{bcd} have proved a similar average theorem for Conjecture \ref{koblitz}:  for any $\ve > 0$, if $\min \{ A(x), B(x) \} \geq x^{1+\ve}$ then
\begin{equation} \label{koblitzthm}
\frac{1}{|\mc{C}(x)|} \sum_{{\begin{substack} {E \in \mc{C}(x)} \end{substack}}} \pi_{E, \text{prime}}(x) \sim C_{\text{prime}} \frac{x}{(\log x)^2}, \quad \textrm{as } x \rightarrow \infty,
\end{equation}
where $C_{\text{prime}}$ is a specific constant.

Finally, Banks and Shparlinski \cite{bs} have proved Conjecture \ref{cyclicity} unconditionally on average:
for any $\ve > 0$, if $x^{2/3+\ve} \leq A(x), B(x) \leq x^{1-\ve}$ then
\begin{equation} \label{cyclicitythm}
\frac{1}{|\mc{C}(x)|} \sum_{{\begin{substack} {E \in \mc{C}(x)} \end{substack}}} \pi_{E, \text{cyclic}}(x) \sim C_{\text{cyclic}} \frac{x}{\log x}, \quad \textrm{as } x \rightarrow \infty,
\end{equation}
where $C_{\text{cyclic}}$ is a specific constant.

In this paper we will prove that each of these average results is consistent with the corresponding conjectured result on the level of the constants.  We will make the notation uniform.
\begin{notation} \label{notation}
Throughout the rest of this paper, let $C_E$ denote any one of the constants $C_{E,r}$, $C_{E,\text{prime}}$, or $C_{E,\text{cyclic}}$, and let $C$ denote the corresponding average constant $C_r$, $C_{\text{prime}}$, or $C_{E,\text{cyclic}}$.
\end{notation}
Our first theorem is conditional upon an affirmative answer to the following question of Serre \cite{serre}.  In its statement, $\mbq(E[p])$ denotes the $p$-th division field of $E$, i.e. the field obtained by adjoining to $\mbq$ the $x$ and $y$-coordinates of a given Weierstrass model of $E$.
\begin{question} \label{serresquestion}
Does there exist an absolute constant $c$ so that, for any prime $p > c$ and any elliptic curve $E$ over $\mbq$ one has
\[
\gal(\mbq(E[p])/\mbq) \simeq GL_2(\mbz/p\mbz)?
\]
\end{question}
We prove
\begin{theorem} \label{maintheoremcond}
Assume that Question \ref{serresquestion} has an affirmative answer.  Then for any positive integer $k$ and any $\ve > 0$, we have
\[
\frac{1}{|\mc{C}|} \sum_{{\begin{substack} {E \in \mc{C}} \end{substack}}} \left| C_E - C \right|^k  \quad \ll_{k,\ve} \quad \max \left\{ \left( \frac{A^\ve \log B}{B} \right)^{\frac{k}{k+1}}, \frac{ \log^{3k + \gamma}(A^3+B^2)}{\sqrt{\min\{A,B\}}} \right\}.
\]
\end{theorem}
Note the following
\begin{corollary}
Provided that Question \ref{serresquestion} has an affirmative answer and that, for some $\ve > 0$,  $A^\ve \ll B \ll e^{A^{\frac{1/2 - \ve}{3+\gamma}}}$, then as $A \rightarrow \infty$, one has
\[
\frac{1}{|\mc{C}|} \sum_{{\begin{substack} {E \in \mc{C}} \end{substack}}} C_{E} \; \longrightarrow \; C
\]
\end{corollary}
Taking $k=2$, we also note the following corollary to Theorem $1.4$ of \cite{dp}, which bounds the average error in the Lang-Trotter conjecture.
\begin{corollary}
Let $\ve > 0$ and $c > 0$ be given and suppose that Question \ref{serresquestion} has an affirmative answer.  Then, provided that $A, B > x^{2+\ve}$ and that
\[
\max \left\{ \left( \frac{A^\ve \log B}{B} \right)^{\frac{2}{3}}, \frac{ \log^{6 + \gamma}(A^3+B^2)}{\sqrt{\min\{A,B\}}} \right\} \; \ll \; \frac{1}{(\log x)^{c-2}},
\]
one has
\[
\frac{1}{|\mc{C}|} \sum_{{\begin{substack} {E \in \mc{C}} \end{substack}}} \left| \pi_{E,r}(x) - C_{E,r}\frac{\sqrt{x}}{\log x} \right|^2 \ll \frac{x}{(\log x)^c}.
\]
\end{corollary}

Unconditionally, we prove a statement about averages over Serre curves (we will review the notion
of a Serre curve in Section \ref{serrecurves}).
\begin{theorem} \label{maintheoremuncond}
Let $k$ be a positive integer.  For each $\ve >0$, we have
\[
\frac{1}{|\mc{C}|} \sum_{{\begin{substack} {E \in \mc{C} \\ E \text{ is a Serre curve }} \end{substack}}} \left| C_{E} - C \right|^k \quad \ll_{k,\ve} \quad A^\ve \cdot \left( \frac{\log B}{B} \right)^{k/(k+1)}.
\]
\end{theorem}
Because of the fact that
\[
\frac{1}{|\mc{C}|} \sum_{{\begin{substack} {E \in \mc{C} \\ E \text{ is a Serre curve }} \end{substack}}} 1 \quad \longrightarrow \quad 1
\]
as $\min\{A,B\} \rightarrow \infty$ (c.f. \cite{jonessc}), Theorem \ref{maintheoremuncond} provides evidence that Theorem \ref{maintheoremcond} should hold unconditionally.

\section{Acknowledgments} \label{acknowledgments}

I wish to thank W. Duke for bringing this problem to my attention and also C. David and A. Granville for comments on an earlier version.

\section{The constants} \label{theconstants}

In this section we will describe precisely the constants occuring in the conjectures under consideration, as well as the corresponding average constants.  Their description involves the division fields of $E$, whose notation we now fix.  
\begin{notation} \label{divisionfieldnotation}
For each positive integer $n$, denote by $\mbq(E[n])$ the $n$-th division field of $E$, obtained by adjoining to $\mbq$ the $x$ and $y$-coordinates of the $n$-torsion of $E$, and by
\[
G_n(E) := \textrm{Gal}\,(\mbq(E[n])/\mbq)
\]
the associated Galois group.  Since $E[n]$ is a free $\mbz/n\mbz$-module of rank $2$, we may (by fixing a $\mbz/n\mbz$-basis) view $G_n(E)$ as a subgroup of $GL_2(\mbz/n\mbz)$.  
\end{notation}
We will distinguish between the case where $E$ has complex multiplication (CM) and the case where $E$ does not (non-CM).  Since
almost all elliptic curves are non-CM \cite{duke}, our only interest in the CM case is to obtain
upper bounds for $C_E$.

\subsection{The CM case} \label{thecmcase}

Suppose that $E$ has complex multiplication by an order $\mc{O}$ in an imaginary quadratic field $K$.  Let $w$ be the number of roots of unity in $K$ and $\gD_K$ the discriminant of $K$.

Then, as computed in \cite[pp. 87--88]{lt}, we have
\begin{equation} \label{ltcmconstant}
C_{E,r} = 
\frac{w}{2\pi} \cdot F_4(r,K) \cdot \prod_{ {\begin{substack} {\ell \neq 2 \\ \ell \mid \gD_K} \end{substack}}} \frac{\ell}{\ell - 1} \cdot \prod_{ {\begin{substack} {\ell \neq 2 \\ \ell \nmid \gD_K \\ \ell \mid r} \end{substack}}} \frac{\ell}{\ell - \left( \frac{\gD_K}{\ell} \right)} \cdot \prod_{ {\begin{substack} {\ell \neq 2 \\ \ell \nmid \gD_K \\ \ell \nmid r} \end{substack}}} \left( 1 - \frac{\left( \frac{\gD_K}{\ell} \right)}{\left(\ell - 1\right)\left(\ell - \left( \frac{\gD_K}{\ell} \right) \right)} \right).
\end{equation}
The factor $F_4(r,K)$ is not relevant to our discussion, and we mention only that, uniformly in $r$ and $K$, one has $| F_4(r,K) | \leq 4$.

To write the Koblitz constant \cite{koblitz}, we define, for any positive integer $n$, the subset
\begin{equation} \label{koblitznotation}
\Phi_n := \{ g \in GL_2(\mbz/n\mbz) : \det(1-g) \in (\mbz/n\mbz)^* \}.
\end{equation}
Then we have
\[
C_{E,\text{prime}} \; = \; \prod_{\ell} \frac{ | G_\ell(E) \cap \Phi_\ell | / | G_\ell(E) |}{1 - 1/\ell} .
\]
Note that (c.f. \cite[Proposition $3.10$]{cojocaruCM}) one has
\[
C_{E,\text{prime}} = \prod_{\ell \mid 6\gD_E} \frac{| G_\ell(E) \cap \Phi_\ell | / | G_\ell(E) |}{1 - 1/\ell} \cdot \prod_{\ell \nmid 6\gD_E} \left( 1 -  \chi(\ell) \frac{\ell^2-\ell-1}{(\ell - \chi(\ell))(\ell-1)^2} \right),
\]
where $\chi$ is the quadratic character corresponding to the quadratic field $K$.

Finally, the cyclicity constant has the same definition on the CM and non-CM case, namely
\begin{equation} \label{cyclicitycmconstant}
C_{E,\text{cyclic}} = \sum_{{\begin{substack} {n \geq 1 } \end{substack}}} \frac{\mu(n)}{[ \mbq(E[n]) : \mbq ]}.
\end{equation}

\subsection{The non-CM case and the average constants}

In each non-CM case, we will write constant $C_E$ in the form
\[
C_E = f(m_E,G_{m_E}(E)) \cdot \prod_{\ell \nmid m_E} f(\ell,GL_2(\mbz/\ell\mbz)),
\]
where $f(n,G)$ is some function of the level $n$ and the subgroup $G \leq GL_2(\mbz/n\mbz)$, and where $m_E$ is a positive integer depending on the torsion representation attached to $E$.  We proceed to describe $m_E$.

Another way to phrase Notation \ref{divisionfieldnotation} is to say that there is a group homomorphism
\[
\varphi_{E,n} : G_\mbq \rightarrow GL_2(\mbz/n\mbz),
\]
defined by letting the absolute Galois group $G_\mbq := \gal (\ol{\mbq}/\mbq)$ act on the $n$-torsion points of $E$, and we are denoting the image of $\varphi_{E,n}$ by $G_n(E)$.  Taking the inverse limit of the $\varphi_{E,n}$ over positive integers $n$ (ordered by divisibility), one obtains a continuous group homomorphism
\[
\varphi_E : G_\mbq \rightarrow GL_2(\hat{\mbz}).
\]
(Here $\hat{\mbz} := \lim_{\leftarrow} \mbz/n\mbz = \prod_p \mbz_p$.)
Serre \cite{serre} showed that, when $E$ has no complex multiplication, the image of this representation is open, i.e. that
\[
[GL_2(\hat{\mbz}) : \varphi_E(G_\mbq)] < \infty.
\]
Equivalently, there is some positive integer level $m_E$ so that, if
\[
\pi : GL_2(\hat{\mbz}) \rightarrow GL_2(\mbz/m_E\mbz)
\]
is the natural projection, one has
\begin{equation} \label{dE}
\varphi_E(G_\mbq) = \pi^{-1}(G_{m_E}(E)).
\end{equation}
For a non-CM curve $E$ over $\mbq$, let us denote by $m_E$ the smallest positive integer such that the above equation holds.  In particular, $m_E$ has the property that, for $m_1$ dividing $m_E$ and $m_2$ coprime to $m_E$ one has
\begin{equation} \label{splits}
G_{m_1m_2}(E) \simeq G_{m_1}(E) \times GL_2(\mbz/m_2\mbz).
\end{equation}

In order to write the Lang-Trotter constant $C_{E,r}$, we follow the notation in \cite{lt}:  for $G \subseteq GL_2(\mbz/n\mbz)$ any subgroup, let
\[
G_r := \{ g \in G \, : \, \tr g \equiv r \mod n \}.
\]
Then,
\[
\begin{split}
C_{E,r} =& \frac{2}{\pi} \cdot \frac{m_E |G_{m_E}(E)_r|}{|G_{m_E}(E)|} \cdot \prod_{{\begin{substack} { \ell \nmid m_E} \end{substack}}} \frac{\ell |GL_2(\mbz/\ell\mbz)_r|}{|GL_2(\mbz/\ell \mbz)|} \\
=& \frac{2}{\pi} \cdot \frac{m_E |G_{m_E}(E)_r|}{|G_{m_E}(E)|} \cdot \prod_{{\begin{substack} {\ell \mid r \\ \ell \nmid m_E} \end{substack}}} \left( 1 + \frac{1}{\ell^2-1} \right) \cdot \prod_{{\begin{substack} {\ell \nmid r \\ \ell \nmid m_E} \end{substack}}} \left( 1 -  \frac{1}{(\ell-1)(\ell^2-1)} \right),
\end{split}
\]
where $m_E$ is as in \eqref{dE}.  On the other hand, the average constant in \eqref{langtrotterthm} is
\[
C_r = \frac{2}{\pi} \cdot \prod_{{\begin{substack} {\ell \mid r } \end{substack}}} \left( 1 + \frac{1}{\ell^2-1} \right) \cdot \prod_{{\begin{substack} {\ell \nmid r} \end{substack}}} \left( 1 -  \frac{1}{(\ell-1)(\ell^2-1)} \right).
\]

The Koblitz constant (as refined by Zywina in \cite{zywina}) is given by
\[
\begin{split}
C_{E, \text{prime}} =& \frac{|G_{m_E}(E) \cap \Phi_{m_E}| / | G_{m_E}(E) |}{\prod_{\ell \mid m_E} (1 - 1/\ell)} \cdot \prod_{\ell \nmid m_E} \frac{|GL_2(\mbz/\ell\mbz) \cap \Phi_{\ell}|/|GL_2(\mbz/\ell \mbz)|}{(1 - 1/\ell)} \\
=& \frac{ |G_{m_E}(E) \cap \Phi_{m_E}|/|G_{m_E}(E)|}{\prod_{\ell \mid m_E} (1 - 1/\ell)} \cdot \prod_{\ell \nmid m_E} \left( 1 - \frac{\ell^2-\ell-1}{(\ell-1)^3(\ell+1)} \right).
\end{split}
\]
In this case the average constant in \eqref{koblitzthm} is given by
\[
C_{\text{prime}} = \prod_{\ell} \left( 1 - \frac{\ell^2-\ell-1}{(\ell-1)^3(\ell+1)} \right)
\]

Finally, the cyclicity constant is given by
\[
C_{E, \text{cyclic}} = \sum_{{\begin{substack} {k \geq 1 } \end{substack}}} \frac{\mu(k)}{|G_k(E)|} = \left( \sum_{{\begin{substack} {k \mid m_E} \end{substack}}} \frac{\mu(k)}{|G_k(E)|} \right) \cdot \prod_{\ell \nmid m_E} \left( 1 - \frac{1}{\ell (\ell-1)^2(\ell+1)} \right),
\]
the second equality coming from \eqref{splits} and the fact that any square-free integer $k$ may be decomposed as $k = k_1 \cdot k_2$, where $k_1 \mid m_E$ and $(k_2,m_E) = 1$.  The average constant in \eqref{cyclicitythm} is
\[
C_{\text{cyclic}} = \prod_{\ell} \left( 1 - \frac{1}{\ell (\ell-1)^2(\ell+1)} \right).
\]

We first note that if any non-CM elliptic curve $E$ were to satisfy $m_E = 1$ (i.e. if $\varphi_E$ were surjective), then we would have $C_E = C$.  However, as observed by Serre, \emph{no} elliptic curve over $\mbq$ has $m_E = 1$.  The main difficulty in proving Theorem \ref{maintheoremcond} is tracking the variation of $m_E$ with $E$.   To prove the theorem, we will focus on a density one subset of curves $E$ for which $m_E$ is essentially equal to the square-free part of the discriminant of $E$.  These curves are called \emph{Serre curves} and will be discussed in detail in the next section.

The proof of Theorem \ref{maintheoremcond} will proceed as follows.   We will decompose the sum
\[
\sum_{E \in \mc{C}} \left| C_E - C \right|^k \; = \; \sum_{{\begin{substack} {E \in \mc{C} \\ E \text{ is a Serre curve}} \end{substack}}} \left| C_E - C \right|^k \; + \; \sum_{{\begin{substack} {E \in \mc{C} \\ E \text{ is not a Serre curve}} \end{substack}}} \left| C_E - C \right|^k.
\]
In Section \ref{averageoverserrecurves}, we will show that, for each $\ve > 0$, one has
\begin{equation*} 
\frac{1}{|\mc{C}|} \sum_{{\begin{substack} {E \in \mc{C} \\ E \text{ is a Serre curve}} \end{substack}}} \left| C_E - C \right|^k \quad \ll_{k,\ve} \quad  A^\ve \cdot \left( \frac{\log B}{B} \right)^{k/(k+1)},
\end{equation*}
proving Theorem \ref{maintheoremuncond}.
In Section \ref{averageovernonserrecurves} we will show that, assuming an affirmative answer to
Question \ref{serresquestion}, one has
\[
\frac{1}{|\mc{C}|} \sum_{{\begin{substack} {E \in \mc{C} \\ E \text{ is not a Serre curve}} \end{substack}}} \left|C_E - C \right|^k \quad \ll \quad \frac{ \log^{3k+\gamma}(A^3+B^2)}{\sqrt{\min\{A,B\}}},
\]
concluding the proof of Theorem \ref{maintheoremcond}.

\section{Serre curves} \label{serrecurves}

Serre \cite{serre} observed that although the torsion representation $\varphi_E$ has finite index in $GL_2(\hat{\mbz})$, it is never surjective when the base field is $\mbq$.  For each elliptic curve $E$ over $\mbq$, there is an index two subgroup $H_E \subseteq GL_2(\hat{\mbz})$ with $\varphi_E(G_\mbq) \subseteq H_E$.  (We will presently describe this subgoup.)
\begin{definition}
An elliptic curve $E$ over $\mbq$ is a \emph{Serre curve} if $\varphi_E(G_\mbq) = H_E$.
\end{definition}
In other words, a Serre curve is an elliptic curve whose torsion representation has image which is ``as large as possible.''

We now describe the subgroup $H_E$:  Let $\gD = \gD_E$ denote the discriminant of $E$ and $\gD_{sf}$ its square-free part, i.e. $\gD_{sf}$ is the unique square-free integer so that 
\[
\frac{\gD}{\gD_{sf}} \in \mbq^2.
\]
The subgroup $H_E$ will be the full preimage under the canonical surjection 
\[
\pi : GL_2(\hat{\mbz}) \rightarrow GL_2(\mbz/M\mbz)
\]
of a particular index two subgroup of $GL_2(\mbz/M\mbz)$ for a certain level $M$.  If $\gD_{sf} = 1$, then every element of $G_2(E)$ must be an even permutation, where we embed $G_2(E)$ into the symmetric group $S_3$ by representing it on the non-identity 2-torsion points.  In this case $M = 2$ and we take $H_E$ to be $\pi^{-1}(A_3)$, where $A_3$ is the alternating group.  Otherwise, the quadratic field $\mbq(\sqrt{\gD})$, being a non-trivial abelian extension of $\mbq$, is contained in a cyclotomic extension.  Let $D_E$ be the smallest positive integer for which 
\[
\mbq(\sqrt{\gD}) \subseteq \mbq(\zeta_{D_E}).
\]
In fact,
\[
D_E = \begin{cases} |\gD_{sf}(E)| & \text{ if } \gD_{sf}(E) \equiv 1 \mod 4 \\
4|\gD_{sf}(E)| & \text{ otherwise } \end{cases},
\]
where $\gD_{sf}(E)$ is the squarefree part of the discriminant of $E$.  Then we define
\begin{equation} \label{defofme}
M_E = \begin{cases} 2|\gD_{sf}(E)| & \text{ if } \gD_{sf}(E) \equiv 1 \mod 4 \\
4|\gD_{sf}(E)| & \text{ otherwise } \end{cases}
\end{equation}
to be the least common multiple of $2$ and $D_E$, so that $\mbq(E[M_E])$ is the compositum of $\mbq(E[2])$ and $\mbq(E[D_E])$.  Since 
\[
\mbq(\sqrt{\gD}) \subseteq \mbq(E[2]) \cap \mbq(E[D_E]),
\]
the corresponding Galois group $G_{M_E}(E) := \textrm{Gal}(\mbq(E[M_E])/\mbq)$ must be a proper subgroup of $GL_2(\mbz/M_E\mbz)$.  In particular, the character on $GL_2(\mbz/M_E\mbz)$ which describes the action on $\sqrt{\gD}$ of an element considering the tower of fields 
\[
\mbq(\sqrt{\gD}) \subseteq \mbq(\zeta_{D_E}) \subseteq \mbq(E[M_E])
\]
is given by 
\[
\gs : \sqrt{\gD} \mapsto \left( \frac{\gD_{sf}(E)}{\det \gs} \right) \sqrt{\gD}.
\]
On the other hand, the inclusion
\[
\mbq(\sqrt{\gD}) \subseteq \mbq(E[2]) \subseteq \mbq(E[M_E])
\]
demands that for each $\gs \in G_{M_E}(E)$, 
\[
\gs : \sqrt{\gD} \mapsto \ve(\gs) \sqrt{\gD},
\]
where $\ve$ denotes the projection $GL_2(\mbz/M_E\mbz) \rightarrow GL_2(\mbz/2\mbz)$ followed by the signature on $GL_2(\mbz/2\mbz) \simeq S_3$.  Thus we see that, with the notation as defined, 
\begin{equation} \label{obstruction}
G_{M_E}(E) \subseteq \ker \left(\ve (\cdot) \left( \frac{\gD_{sf}(E)}{\det(\cdot)} \right) \right) \subseteq GL_2(\mbz/M_E\mbz).
\end{equation}
In this case we therefore make the definition 
\[
H_E = \pi^{-1}\left( \ker \left(\ve (\cdot) \left( \frac{\gD_{sf}(E)}{\det(\cdot)} \right) \right) \right).
\]
Note that $M_E$ always divides $m_E$, and if we interpret $\left( \frac{\gD_{sf}(E)}{\det(\cdot)} \right)$ to be the trivial character in case $\gD_{sf} = 1$, we have
\begin{equation} \label{meequalsME}
E \text{ is a Serre curve } \; \Longleftrightarrow \; m_E = M_E \text{ and } G_{M_E}(E) = \ker \left(\ve (\cdot) \left( \frac{\gD_{sf}(E)}{\det(\cdot)} \right) \right).
\end{equation}
One shows easily that, for $d$ a proper divisor of  $M_E$ and $\pi$ denoting the natural projection $GL_2(\mbz/M_E\mbz) \rightarrow GL_2(\mbz/d\mbz)$, one has
\[
\pi \left( \ker \left(  \ve \cdot \left( \frac{\gD_{sf}(E)}{\det (\cdot)} \right) \right) \right) = GL_2(\mbz/d\mbz).
\]
Thus in particular, when $E$ is a Serre curve and $d \mid M_E$, one has
\begin{equation} \label{meminimalexceptional}
G_d(E) = \begin{cases}
		\ker \left(  \ve \cdot \left( \frac{\gD_{sf}(E)}{\det (\cdot)} \right) \right) & \text{ if } d = M_E \\
		GL_2(\mbz/d\mbz) & \text{ otherwise.}
		\end{cases}
\end{equation}

\section{The average over Serre curves} \label{averageoverserrecurves}

We will now show that for each $\ve > 0$, one has
\begin{equation} \label{overserrecurves}
\frac{1}{|\mc{C}|} \sum_{{\begin{substack} {E \in \mc{C} \\ E \text{ is a Serre curve}} \end{substack}}} \left| C_E - C \right|^k \quad \ll_{k,\ve} \quad A^\ve \cdot \left( \frac{\log B}{B} \right)^{k/(k+1)}.
\end{equation}

\subsection{The constants associated to Serre curves}  \label{serrecurveconstants}

In this section, we will explicitly compute the constants $C_{E,r}$, $C_{E, \text{prime}}$ and $C_{E, \text{cyclic}}$ for $E$ a Serre curve.  For the Lang-Trotter constant $C_{E,r}$, we must fix some notation.  First define the exponent $k \in \{ 1,2,3 \}$ and the odd integer $W$ by
\begin{equation*} \label{defofW}
W := \frac{\gD_{sf}}{(\gD_{sf},2)} \quad\quad \text{ and } \quad\quad
									k := 	\begin{cases}
										1 & \text{ if } \gD_{sf} \equiv 1 \mod 4 \\
										2 & \text{ if } \gD_{sf} \equiv 3 \mod 4 \\
										3 & \text{ if } \gD_{sf} \equiv 2 \mod 4.
										\end{cases}
\end{equation*}
In other words, we have
\begin{equation} \label{2kW}
M_E =: 2^k \cdot W,
\end{equation}
where $M_E$ is as in \eqref{defofme}.  When $2^{k-1}$ divides $r$, we further define the symbol $\gd(\gD_{sf},r) \in \{ \pm 1 \}$ by
\begin{equation} \label{defofgd}
\gd(\gD_{sf},r) := (-1)^{\om \left( \frac{W}{(W,r)} \right) + \frac{W+1}{2} + \frac{r}{2^{k-1}}} \cdot \chi_4 \left( -\frac{\gD_{sf}}{2} \right),
\end{equation}
where we make the convention that $\chi_4(x) = 1$ if $x \notin \mbz$.

\begin{proposition} \label{serreconstantprop}
Suppose that $E$ is an elliptic curve over $\mbq$ which is a Serre curve.  Then
\begin{equation} \label{cer}
C_{E,r} = \begin{cases} 		C_r \left( 1 + \gd(\gD_{sf},r) \cdot \frac{ M_E \cdot 2^{k-1} \cdot \varphi( (W,r) )}{ | 							GL_2(\mbz/M_E\mbz)_r |}  \right) & \text{ if } 2^{k-1} \mid r \\
					C_r & \text{ otherwise,}
\end{cases}
\end{equation}
\begin{equation} \label{ceprime}
C_{E, \text{prime}} = \begin{cases}  C_{\text{prime}} \left( 1 + \prod_{p \mid \gD_{sf}} \frac{1}{p^3-2p^2-p+3} \right) & \text{ if } \gD_{sf} \equiv 1 \mod 4 \\
C_{\text{prime}} & \text{ otherwise}
\end{cases},
\end{equation}
and
\begin{equation} \label{cecyclic}
C_{E, \text{cyclic}} = \begin{cases}
				C_{\text{cyclic}} \left( 1 + \frac{\mu(M_E)}{\prod_{\ell \mid M_E} \left( |GL_2(\mbz/\ell\mbz)| - 1 \right)} \right) & if \gD_{sf} \equiv 1 \mod 4 \\
				C_{\text{cyclic}} & \text{ otherwise} \end{cases}.
\end{equation}
\end{proposition}
The proof of the proposition will require the use of some technical lemmas.
We now describe the set-up of the first of these lemmas.  Let $M$ be any positive integer and
recall the isomorphism of the Chinese remainder theorem:
\begin{equation} \label{CRT}
GL_2(\mbz/M\mbz) \simeq \prod_{p^k \mid\mid M} GL_2(\mbz/p^k\mbz), \quad \quad x \mapsto (x_{p^k})
\end{equation}
Suppose that $X_M \subseteq GL_2(\mbz/M\mbz)$ is any subset which, under \eqref{CRT}, satisfies
\[
X_M \simeq \prod_{p^k \mid\mid M} X_{p^k},
\]
where $X_{p^k}$ denotes the projection of $X_M$ onto the $p^k$-th factor.  Suppose further that,
for each prime $p$ dividing $M$ we have a group homomorphism
\[
\psi_{p^k} : GL_2(\mbz/p^k\mbz) \longrightarrow \{ \pm 1 \},
\]
and write
\[
\psi_M : GL_2(\mbz/M\mbz) \longrightarrow \{ \pm 1 \}, \quad \quad \psi_M(x) := \prod_{p^k \mid\mid M}\psi_{p^k}(x_{p^k}).
\]
\begin{lemma} \label{psi}
With notation as just outlined, we have
\[
|\psi_M^{-1}(\pm 1) \cap X_M| \; = \; \frac{1}{2} \left( | X_M |  \; \pm \; \prod_{p^k \mid\mid M} \left( | \psi_{p^k}^{-1}(1) \cap X_{p^k} | - | \psi_{p^k}^{-1}(-1) \cap X_{p^k} | \right) \right).
\]
\end{lemma}
\begin{proof}
We begin by noting that
\begin{equation*} \label{psiformula}
\begin{split}
\left| \psi_M^{-1}(\pm 1) \cap X_M \right| \quad =& \quad \sum_{{\begin{substack} {(s_p)_{p \mid M} \\ \prod s_p = \pm 1} \end{substack}}} \prod_{p^k \mid\mid M} \left| \psi_{p^k}^{-1}(s_p) \cap X_{p^k} \right| \\
=& \sum_{{\begin{substack} {(s_p)_{p \mid M} \\ \prod s_p = \pm 1} \end{substack}}} \prod_{p^k \mid\mid M} \left( F_{1}(p^k) + s_pF_{-1}(p^k) \right),
\end{split}
\end{equation*}
where
\[
F_{1}(p^k) := \frac{1}{2} |X_{p^k}| \quad \text{ and } \quad F_{-1}(p^k) := \frac{1}{2} \left( | \psi_{p^k}^{-1}(1) \cap X_{p^k} | - | \psi_{p^k}^{-1}(-1) \cap X_{p^k} | \right).
\]
Here our notation is meant to indicate that the sum runs over all $\om(M)$-tuples $(s_p)_{p \mid M}$ of $\pm 1$'s which satisfy
$\prod_{p \mid M} s_p = \pm 1$.  Expanding the product and reversing summation, we obtain
\begin{equation} \label{reversedsum}
\left| \psi_M^{-1}(\pm 1) \cap X_M \right| \quad = \quad \sum_{{\begin{substack} {(t_p)_{p \mid M}} \end{substack}}} \prod_{p^k \mid\mid M} F_{t_p}(p^k) \left( \sum_{{\begin{substack} {(s_p)_{p \mid M} \\ \prod s_p = \pm 1} \end{substack}}} \left( \prod_{{\begin{substack} {p \mid M \\ t_p = -1} \end{substack}}} s_p \right) \right),
\end{equation}
where now $(t_p)$ runs over \emph{all} $\om(M)$-tuples of $\pm 1$'s.  Now we show that, for all tuples $(t_p)$ except $(t_p) \in \{ (1,1,\dots,1), (-1,-1,\dots,-1)\}$, the innermost sum is equal to zero.  For suppose that
\[
\{ p : p \mid M, t_p = 1 \} \; \neq \; \emptyset \; \neq \; \{ p : p \mid M, t_p = -1 \},
\]
and fix a prime $p_1$ with $t_{p_1} = 1$ and a prime $p_2$ with $t_{p_2} = -1$.  For a tuple $(s_p)$, define its dual $(\hat{s}_p)$ by
\[
\hat{s}_{p_1} = -s_{p_1}, \quad \hat{s}_{p_2} = -s_{p_2}, \quad \text{ and } \quad \hat{s}_{p} = s_{p} \quad (p \notin \{ p_1, p_2 \}).
\]
Noting that 
\[
\prod_{p \mid M } s_p = \prod_{p \mid M} \hat{s}_p  \quad \text{ and } \quad \prod_{{\begin{substack} {p \mid M\\ t_p = -1} \end{substack}}} s_p + \prod_{{\begin{substack} {p \mid M \\ t_p = -1} \end{substack}}} \hat{s}_p = 0,
\]
we see that, except when $(t_p) \in \{ (1,1,\dots,1), (-1,-1,\dots,-1)\}$, the innermost sum in \eqref{reversedsum} vanishes.  Thus,
\[
\left| \psi_M^{-1}(\pm 1) \cap X_M \right| \quad = \quad \frac{1}{2} \left( \prod_{p^k \mid\mid M} F_{1}(p^k) \pm \prod_{p^k \mid\mid M} F_{-1}(p^k) \right),
\]
proving the lemma.
\end{proof}

\noindent \emph{Proof of \eqref{cer}.}
If $E$ is a non-CM curve then
\[
\frac{C_{E,r}}{C_r} \quad = \quad  \frac{m_E | G_{m_E}(E)_r |}{ | G_{m_E}(E) |} \cdot \prod_{\ell \mid m_E} \frac{| GL_2(\mbz/\ell\mbz) |}{ \ell | GL_2(\mbz/\ell\mbz)_r |}.
\]
Thus, when $E$ is a Serre curve, we use \eqref{meequalsME} and \eqref{2kW} to write
\begin{equation*} 
\frac{C_{E,r}}{C_r} \quad = \quad \frac{2^{k} W}{2W} \cdot \frac{2 \left| \left( \ker\left( \ve \cdot \left( \frac{\gD_{sf}}{\det(\cdot)} \right) \right) \right)_r \right|}{ | GL_2(\mbz/M_E\mbz)_r |} \cdot \frac{| GL_2(\mbz/2^k\mbz)_r |}{ | GL_2(\mbz/2\mbz)_r | } \cdot \frac{ | GL_2(\mbz/2\mbz) | }{ | GL_2(\mbz/2^k\mbz) |}.
\end{equation*}
Now we use the following corollary of Lemma \ref{chi} below.
\begin{corollary} \label{sizeofgl22k}
For $k \in \{ 1, 2, 3 \}$, we have
\[
| GL_2(\mbz/2^k\mbz)_r | = \begin{cases} 2^{3k-1} & \text{ if $r$ is even } \\
								2^{3k-2} & \text{ if $r$ is odd.}
								\end{cases}
\]
\end{corollary}
The corollary, together with $| GL_2(\mbz/2^k\mbz) | = 3\cdot 2^{4k-3}$, implies that
\begin{equation*}
\frac{C_{E,r}}{C_r} \quad = \quad \frac{2 \left| \left( \ker\left( \ve \cdot \left( \frac{\gD_{sf}}{\det(\cdot)} \right) \right) \right)_r \right|}{ | GL_2(\mbz/M_E\mbz)_r |}.
\end{equation*}
To evaluate $\left| \left( \ker\left( \ve \cdot \left( \frac{\gD_{sf}}{\det(\cdot)} \right) \right) \right)_r \right|$, we will apply Lemma \ref{psi} with $M = M_E$ and
\begin{equation} \label{defofpsi}
\psi_{p^k}(\gs) := \begin{cases}
               	\left( \frac{\det(\gs)}{p} \right) & \text{ if $p$ is odd} \\
		\ve(\gs) & \text{ if $p^k = 2$ and } \gD_{sf} \equiv 1 \mod 4 \\
\chi_4(\det \gs) \ve(\gs) & \text{ if $p^k = 4$ and } \gD_{sf} \equiv 3 \mod 4 \\
\chi_8(\det \gs) \ve(\gs) & \text{ if $p^k = 8$ and } \gD_{sf} \equiv 2 \mod 8 \\
\chi_4(\det \gs) \chi_8(\det \gs) \ve(\gs) & \text{ if $p^k = 8$ and } \gD_{sf} \equiv 6 \mod 8.
              \end{cases}
\end{equation}
Note that we then have $\ve \cdot \left( \frac{\gD_{sf}}{\det(\cdot)} \right) = \prod_{p^k \mid\mid M_E} \psi_{p^k}(\cdot)$.  Thus, by Lemma \ref{psi}, we have
\begin{equation} \label{cerbycr}
\frac{C_{E,r}}{C_r} \quad = \quad 1 + \frac{\prod_{p^k \mid\mid M_E}\left( |\psi_{p^k}^{-1}(1)_r| - |\psi_{p^k}^{-1}(-1)_r \right)}{ | GL_2(\mbz/M_E\mbz)_r |}.
\end{equation}

\begin{lemma} \label{podderror}
For odd primes $p$, one has
\[
\psi_{p^k}(1)_r - \psi_{p^k}(-1)_r = \begin{cases}
                                      	\left( \frac{-1}{p} \right) p(p-1) & \text{ if } p \mid r \\
					- \left( \frac{-1}{p} \right) p & \text{ if } p \nmid r.
                                     \end{cases}
\]
\end{lemma}
\noindent \emph{Proof of Lemma \ref{podderror}.}
For $r$ and $d$ modulo $p$, it is straightforward to show that 
\[
\left| \left\{ g \in GL_2(\mbz/p\mbz) : \tr g = r, \; \det g = d \right\} \right| \; = \; p\left(p+\left(\frac{r^2-4d}{p}\right)\right).
\]
Thus, partitioning $\psi_p^{-1}(\pm 1)_r$ by determinant shows that 
\begin{equation}
\label{psipformula}
|\psi_p^{-1}(\pm 1)_r| = \sum_{\left(\frac{d}{p}\right) = \pm 1} p\left(p+\left(\frac{r^2-4d}{p}\right)\right).
\end{equation}
To evaluate this sum we use the following corollary of \cite[Lemma $6$]{jonestf}.
\begin{sublemma}
For $r$ nonzero modulo $p$,
\[
\sum_{\left(\frac{d}{p}\right) = \pm 1} \left(\frac{r^2-4d}{p}\right) \; = \; -\frac{1 \pm \left(\frac{-1}{p}\right)}{2},
\]
while if $r \equiv 0 \mod p$ then
\[
\sum_{\left(\frac{d}{p}\right) = \pm 1} \left(\frac{r^2-4d}{p}\right) \; = \; \pm \left(\frac{-1}{p}\right) \frac{p-1}{2}.
\]
\end{sublemma}
Inserting this into \eqref{psipformula}, we finish the proof of Lemma \ref{podderror}. \hfill $\Box$ \\

For $p = 2$, we use the following lemma, whose proof is a straightforward calculation which we omit.
\begin{lemma} \label{chi}
If $\gD_{sf} \equiv 1 \mod 4$ then
\[
|\psi_2^{-1}(1)_0| = 1, \; |\psi_2^{-1}(-1)_0| = 3, \; |\psi_2^{-1}(1)_1| = 2, \; \text{ and } \, |\psi_2^{-1}(-1)_1| = 0.
\]
If $\gD_{sf} \equiv 3 \mod 4$ then
\[
|\psi_4^{-1}(1)_0| = |\psi_4^{-1}(-1)_2| = 12, \; |\psi_4^{-1}(-1)_0| = |\psi_4^{-1}(1)_2| = 20,
\]
and for any odd $r$ modulo $4$,
\[
|\psi_4^{-1}(\pm 1)_r| = 8.
\]
If $\gD_{sf} \equiv 2 \mod 4$ then
\[
|\psi_8^{-1}(\pm 1)_r| = \begin{cases}
16 \cdot 8 & \text{ if } r \equiv 2 \mod 4 \\ 
16 \cdot 4 & \text{ if } r \textrm{ is odd. } 
\end{cases},
\]
while
\[
|\psi_8^{-1}(1)_0| = |\psi_8^{-1}(-1)_4| = \begin{cases}
16 \cdot 9 & \textrm{ if } \gD_{sf} \equiv 2 \mod 8 \\
16 \cdot 7 & \textrm{ if } \gD_{sf} \equiv 6 \mod 8
\end{cases}
\]
and 
\[
|\psi_8^{-1}(-1)_0| = |\psi_8^{-1}(1)_4| = \begin{cases}
16 \cdot 7 & \textrm{ if } \gD_{sf} \equiv 2 \mod 8 \\
16 \cdot 9 & \textrm{ if } \gD_{sf} \equiv 6 \mod 8.
\end{cases}
\]
\end{lemma}
\begin{corollary} \label{chicorollary}
For $\psi_{2^k}$ as in \eqref{defofpsi}, we have
\[
\psi_{2^k}(1)_r - \psi_{2^k}(-1)_r = \begin{cases}
                                      	- (-1)^{r/2^{k-1}} \chi_4(-\gD_{sf}/2) \cdot 2^{2k-1} & \text{ if } 2^{k-1} \mid r \\
					0 & \text{ otherwise,}
                                     \end{cases}
\]
where here we use the convention that $\chi_4(x) = 1$ if $x$ is not an integer.
\end{corollary}
Inserting the results of Corollary \ref{chicorollary} and Lemma \ref{podderror} into \eqref{cerbycr}, we finish the proof of \eqref{cer}. \hfill $\Box$

Having proved \eqref{cer}, we now proceed to \\

\noindent \emph{Proof of \eqref{ceprime}.}
This computation may also be found in \cite{zywina}.  For any non-CM elliptic curve $E$, we have
\[
\frac{C_{E,\text{prime}}}{C_{\text{prime}}} = \frac{| G_{m_E}(E) \cap \Phi_{m_E} |}{|G_{m_E}(E)|} \cdot \prod_{\ell \mid m_E} \left( \frac{ | GL_2(\mbz/\ell\mbz) | }{ | \Phi_\ell | } \right).
\]
If $E$ is a Serre curve, then we have
\begin{equation} \label{ceprimebycprime}
\frac{C_{E,\text{prime}}}{C_{\text{prime}}} =  \frac{2 | \psi_{M_E}^{-1}(1) \cap \Phi_{M_E} |}{ | \Phi_{M_E} |}.
\end{equation}
Applying Lemma \ref{psi}, we find that
\[
 | \psi_{M_E}^{-1}(1) \cap \Phi_{M_E} | = \frac{1}{2} \left( |\Phi_{M_E} | + \prod_{p^k \mid\mid M_E} \left( | \psi_{p^k}^{-1}(1) \cap \Phi_{p^k} | - |  \psi_{p^k}^{-1}(-1) \cap \Phi_{p^k} | \right) \right).
\]
\begin{lemma} \label{psipkoblitz}
For $p$ odd, one has
\[
|\Phi_p | = p \left( p^3 - 2p^2 - p + 3 \right)
\]
and
\[
| \psi_{p}^{-1}(1) \cap \Phi_{p} | - |  \psi_{p}^{-1}(-1) \cap \Phi_{p} | = p.
\]
\end{lemma}
\noindent \emph{Proof of Lemma \ref{psipkoblitz}.}
The lemma follows immediately from
\[
|\psi_{p}^{-1}(\pm 1) \cap \Phi_{p} | = \frac{1}{2} \cdot p \left( p^3 - 2p^2 - p + 3 \pm 1 \right).
\]
\hfill $\Box$ \\
The following lemma is a straightforward calculation using the fact that
\[
\Phi_{2^k} = \ve^{-1}(1).
\]
\begin{lemma} \label{psi2koblitz}
One has
\[
| \psi_{2^k}^{-1}(1) \cap \Phi_{2^k} | - |  \psi_{2^k}^{-1}(-1) \cap \Phi_{2^k} | = \begin{cases}
2 & \text{ if } k=1 \\
0 & \text{ if } k \in \{2,3\}.
\end{cases}
\]
\end{lemma}
Inserting the results of Lemmas \ref{psipkoblitz} and \ref{psi2koblitz} into \eqref{ceprimebycprime}, we finish the proof of \eqref{ceprime}. \hfill $\Box$ \\

\noindent \emph{Proof of \eqref{cecyclic}.}
We have
\[
\frac{C_{E,\text{cyclic}}}{C_{\text{cyclic}}} = \frac{\sum_{k \mid m_E} \frac{\mu(k)}{ | G_k(E) | } }{\prod_{\ell \mid m_E} \left( 1 - \frac{1}{ | GL_2(\mbz/\ell\mbz) | } \right) }.
\]
If $E$ is a Serre curve then $m_E = M_E$.  Note that $M_E$ is square-free if and only if $\gD_{sf}(E) \equiv 1 \mod 4$.  Thus, if $E$ is a Serre curve, we deduce from \eqref{meminimalexceptional} that
\[
\sum_{k \mid m_E} \frac{\mu(k)}{ | G_k(E) | }  = \begin{cases}
									\prod_{\ell \mid m_E} \left( 1 - \frac{1}{| GL_2(\mbz/\ell										\mbz) |} \right) + \frac{\mu(m_E)}{|GL_2(\mbz/m_E\mbz)|}   									& \text{ if } \gD_{sf} \equiv 1 \mod 4 \\
									\prod_{\ell \mid m_E} \left( 1 - \frac{1}{| GL_2(\mbz/\ell										\mbz) |} \right) & \text{ otherwise. } 									\end{cases}
\]
This proves \eqref{cecyclic}.
\hfill $\Box$ \\

We have now proved \eqref{cer}, \eqref{ceprime}, and \eqref{cecyclic}, finishing the proof of the Proposition \ref{serreconstantprop}.

\subsection{Averaging the Serre curve constants.} 

Considering Proposition \ref{serreconstantprop}, we see that when $E$ is a Serre curve, $C_E$ has the form
\[
C_E = C \left( 1 + g(\gD_{sf}(E)) \right) \quad \text{ where } \quad g(\gD_{sf}(E)) \ll \frac{1}{ \gD_{sf}(E) }.
\]
Since the discriminant of the curve $Y^2=X^3+aX+b$ is $-16(4a^3+27b^2)$, the result \eqref{overserrecurves} will follow from
\begin{equation} \label{afterZ}
\frac{1}{4AB}\sum_{{\begin{substack} {|a| \leq A \\ |b| \leq B \\ 4a^3+27b^2 \neq 0} \end{substack}}} \frac{1}{|(4a^3+27b^2)_{sf}|^k} \; \ll_\ve \; A^\ve \cdot \left( \frac{\log B}{B} \right)^{k/(k+1)}.
\end{equation}
Let $Z$ be a positive real number to be chosen later.  Since the left hand side is bounded by
\[
\frac{1}{4AB}\sum_{{\begin{substack} {|a| \leq A \\ |b| \leq B \\ 4a^3+27b^2 \neq 0 \\ |(4a^3+27b^2)_{sf}| \leq Z} \end{substack}}} 1 \quad +  \quad \frac{1}{4AB}\sum_{{\begin{substack} {|a| \leq A \\ |b| \leq B \\ |(4a^3+27b^2)_{sf}| > Z } \end{substack}}} \frac{1}{Z^k},
\]
we are reduced to proving the following lemma.
\begin{lemma} \label{squarefreesievelemma}
For any $\ve > 0$, we have
\begin{equation} \label{squarefreesum}
\sum_{{\begin{substack} {|a| \leq A \\ |b| \leq B \\ 4a^3+27b^2 \neq 0 \\ |(4a^3+27b^2)_{sf}| \leq Z} \end{substack}}} 1 \; \ll_\ve \; \log B \cdot Z \cdot A^{1+\ve}.
\end{equation}
\end{lemma}
\begin{proof}
The proof boils down to counting ideals of bounded norm in various quadratic fields.  I would like to thank R. Daileda for helpful discussions regarding this viewpoint.  We wish to count the number of integer pairs $(a,b) \in [-A,A]\times [-B,B]$ which satisfy the equation
\[
4a^3+27b^2 = dy^2,
\]
where $y$ and $d$ are integers with $d \neq 0$ square-free and $|d| \leq Z$.  Re-writing this equation as
\[
x^2 - D y^2 = 12(-a)^3,
\]
where $x = 9b$ and $D = 3d$, we see that the left hand side of \eqref{squarefreesum} is bounded by
\begin{equation} \label{boundone}
\sum_{{\begin{substack} { 2 < |D| \leq 3Z \\ D \text{ square-free} } \end{substack}}} \sum_{1 \leq |a| \leq A} \# \{ (\ga,\ol{\ga}) \in \left( \mc{O}_{\mbq(\sqrt{D})} \right)^2 : \ga \cdot \ol{\ga} = 12a^3, \; \ga + \ol{\ga} \leq 18B \},
\end{equation}
where $\mc{O}_{\mbq(\sqrt{D})}$ denotes the ring of integers of $\mbq(\sqrt{D})$.  We will presently transform this into counting principal ideals rather than elements up to conjugation, but in the real quadratic case we must worry about the presence of an infinite unit group.  Suppose that $(\ga,\ol{\ga}) \in \left( \mc{O}_{\mbq(\sqrt{D})} \right)^2$ is a conjugate pair satisfying
\[
\ga \cdot \ol{\ga} = 12a^3 \quad \text{ and } \quad | \ga + \ol{\ga} | \leq 18B.
\]
Writing $\ga = r + s\sqrt{D}$, we may assume that $r$ and $s$ have the same sign.
Any other $\gb \in  \mc{O}_{\mbq(\sqrt{D})}$ which generates the same ideal as $\ga$ is of the form $\gb = \ga \cdot \ve_D^n$, where $\ve_D$ is a fundamental unit.  One can show that the number of integers $n$ for which
\[
| \ga \cdot \ve_D^n + \ol{\ga} \cdot \ol{\ve_D}^n | \leq 18B
\]
is $\ll \log B$, with an absolute constant.  Thus, \eqref{boundone} is bounded by a constant times
\[
\log B \cdot \sum_{{\begin{substack} {2 <  |D| \leq 3Z \\ D \text{ square-free} } \end{substack}}} \sum_{1\leq a \leq A} \eta_D^{\text{princ}}(12a^3) \; \leq \; \log B \cdot \sum_{{\begin{substack} {2 < |D| \leq 3Z \\ D \text{ square-free} } \end{substack}}} \sum_{1 \leq a \leq A} \eta_D(12a^3),
\]
where $\eta_D(m)$ (resp. $\eta_D^{\text{princ}}(m)$) is the number of integral ideals (resp. the number of \emph{principal} ideals) in the ring $\mc{O}_{\mbq(\sqrt{D})}$ of norm equal to $m$.  We will now show that
\begin{equation} \label{taubound}
\eta_D(m) \leq \tau(m).
\end{equation}
To see this, note that the set of integral ideals $I$ of $\mc{O}_{\mbq(\sqrt{D})}$ of norm $m$ is exactly
\[
\{ \prod_{{\begin{substack} { p^{\ga_p} \mid\mid m \\ p \text{ split} } \end{substack}}} \mathfrak{P}^i \ol{\mathfrak{P}}^{\ga_p-i}\cdot \prod_{{\begin{substack} { p^{\ga_p} \mid\mid m \\ p \text{ inert} } \end{substack}}} (p \mc{O}_{\mbq(\sqrt{D})})^{\ga_p/2} \cdot \prod_{{\begin{substack} { p^{\ga_p} \mid\mid m \\ p \text{ ramified} } \end{substack}}} \mathfrak{P}^{\ga_p} \},
\]
where $\mathfrak{P}$ denotes a prime ideal lying over $p$ and $0 \leq i \leq \ga_p$.  The number of such choices is
\[
\eta_D(m) = 
	\begin{cases}
		\prod_{{\begin{substack} { p^{\ga_p} \mid\mid m \\ p \text{ split} } \end{substack}}} (\ga_p+1) & \text{ if } p \text{ inert } \Rightarrow 2 \mid \ga_p, \\
		0 & \text{ otherwise.} 
	\end{cases}
\]
From this, \eqref{taubound} is immediate.
Noting that $\tau(m) \ll_\ve m^\ve$, one sees that
\[
\log B \cdot \sum_{{\begin{substack} {2 < |D| \leq 3Z \\ D \text{ square-free} } \end{substack}}} \sum_{1 \leq a \leq A} \eta_D(12a^3) \; \ll_\ve \; \log B \cdot Z \cdot A^{1+\ve},
\]
finishing the proof of Lemma \ref{squarefreesievelemma}.
\end{proof}
Finally, using $Z = (B/(A^\ve \log B))^{1/(k+1)}$, \eqref{afterZ} follows, and thus so does \eqref{overserrecurves}.

\section{The average over non-Serre curves} \label{averageovernonserrecurves}
We finally turn to proving that
\begin{equation} \label{overnonserrecurves}
\frac{1}{|\mc{C}|} \sum_{{\begin{substack} {E \in \mc{C} \\ E \text{ is not a Serre curve}} \end{substack}}} \left| C_E - C \right|^k \ll_k \frac{ \log^{3k+\gamma}(A^3+B^2)}{\sqrt{\min\{A,B\}}} .
\end{equation}
In the case of the cyclicity constant, one has
\[
C_{E,\text{cyclic}} \leq 1,
\]
since the constant is a density.  For the other constants, we prove the following.
\begin{lemma} \label{constantbounds}
Fix the integer $r \in \mbz$.
If $E$ over $\mbq$ is an elliptic curve with CM, then
\[
C_{E,r} = O(1) \quad \text{ and } \quad  C_{E, prime} = O(\log(\gD_E)).
\]
If $E$ is a non-CM elliptic curve over $\mbq$, then we have
\[
C_{E, prime} = O(\log m_E).
\]
Assuming an affirmative answer to Question \ref{serresquestion}, we have
\[
C_{E, r} = O(\log^3 m_E).
\]
\end{lemma}
\begin{proof}
For the Lang-Trotter constant in the CM case, we see from \eqref{ltcmconstant} that
\[
C_{E,r} \ll \prod_{{\begin{substack} {\ell \neq 2 \\ \ell \mid \gD_K} \end{substack}}} \left( 1 + \frac{1}{\ell - 1} \right) \cdot \prod_{{\begin{substack} {\ell \neq 2 \\ \ell \nmid \gD_K \\ \ell \mid r} \end{substack}}} \left( 1 +\frac{1}{\ell - 1} \right) \ll \log \gD_K \cdot \log r,
\]
by Merten's theorem.  For fixed $r$ and using the fact (c.f. \cite{silverman2}) that, since $E$ is defined over $\mbq$, $K$ must have class number one, we conclude that $C_{E,r}$ is uniformly bounded. 

For $C_{E,r}$ in the non-CM case, we reason as follows.  According to \cite{jonestc}, we may take $m_E$ to be of the form
\[
m_E = \prod_{{\begin{substack} { p \in \{ 2,3,5 \} \text{ or } \\ G(p) \subsetneq GL_2(\mbz/p\mbz)} \end{substack}}} p^{\ga_p} \cdot \prod_{{\begin{substack} { p \mid \gD_E \\ p \notin \{ 2, 3,5\} \text{ and } \\ G(p) = GL_2(\mbz/p\mbz)} \end{substack}}} p^{\ga_p} =: m_1 \cdot m_2,
\]
where $\ga_p$ are certain exponents, independent of $E$.  If Question \ref{serresquestion} has an affirmative answer then $m_1$ is uniformly bounded, and there must be a non-trivial intersection
\[
\mbq \subsetneq \mbq(E[m_1]) \cap \mbq(E[m_2]) =: F,
\]
whose Galois group we will denote by $H$:
\[
H := \gal(F/\mbq).
\]
The restriction of an automorphism to the subfield $F$ defines group homomorphisms
\[
\chi : G(m_1) \longrightarrow H, \quad \quad \psi : G(m_2) \longrightarrow H,
\]
and the Galois group of the $m_E$-th division field may be identified as
\[
\gal(\mbq(E[m_E)/\mbq) = \{ (\tau_1,\tau_2) \in G(m_1) \times G(m_2) : \chi (\tau_1) = \psi (\tau_2 ) \} = \ker \chi \otimes \psi^{-1}.
\]
Because of basic facts about the group $GL_2(\mbz_p)$ (c.f. \cite{serre2}) one has $G(m_2) = GL_2(\mbz/m_2\mbz)$.  As observed in \cite{serre2}, the common quotient $H$ of $G(m_1)$ and $GL_2(\mbz/m_2\mbz)$ must be abelian, and since the commutator subgroup of $GL_2(\mbz/m_2\mbz)$ is equal to $SL_2(\mbz/m_2\mbz)$, we see that $\psi = \gd \circ \det$ for some homomorphism
\[
\gd : (\mbz/m_2\mbz)^* \longrightarrow H.
\]
Considering the decomposition
\[
G(m_E)_r = \left( \ker \chi \otimes \psi^{-1} \right)_r = \bigsqcup_{h \in H} \chi^{-1}(h)_r \times \det{}^{-1}(\gd^{-1}(h))_r,
\]
we are led to
\begin{lemma} \label{bonnlemma}
Fix any odd prime power $p^{n}$ and integers $r \in \mbz/p^n\mbz$ and $d \in (\mbz/p^n\mbz)^*$.  Then we have
\[
| \{ A \in M_{2 \times 2}(\mbz/p^n\mbz) : \tr(A) = r, \det(A) = d \} | \; \leq \; p^{2n} \left( 1 + \frac{3}{p} \right)
\]
\end{lemma}
\noindent \emph{Proof of Lemma \ref{bonnlemma}.}
In fact, we will evaluate the left-hand side explicitly.  Writing
\[
r^2-4d =: \gD =: p^\gd \cdot \gD',
\]
where $p \nmid \gD'$, we will show that
\begin{equation} \label{exactformula}
| \{ A \in M_{2 \times 2}(\mbz/p^n\mbz) : \tr(A) = r, \det(A) = d \} | \; = \; p^{2n} \left( 1 + \frac{1}{p} + f(p,\gD) \right),
\end{equation}
where
\[
f(p,\gD) := \begin{cases}
             	-p^{-(n+1)/2} & \text{ if } \gd = n \text{ is odd} \\
		-p^{-(n+2)/2} & \text{ if } \gd = n \text{ is even} \\
		-\left( p^{-(\gd+1)/2} + p^{-(\gd+3)/2} \right) & \text{ if } \gd < n, 2 \nmid \gd \\
		\left(\left( \frac{\gD/p^\gd}{p} \right) \left( \gd + 2 - \frac{\gd+1}{p} \right) + \gd - \frac{\gd+1}{p} \right) p^{-(\gd/2+1)} & \text{ if } \gd < n, 2 \mid \gd,
            \end{cases}
\]
where $\left( \frac{\cdot}{p} \right)$ denotes the Legendre symbol.  Lemma \ref{bonnlemma} follows upon observing that
\[
\left| \left(\left( \frac{\gD/p^\gd}{p} \right) \left( \gd + 2 - \frac{\gd+1}{p} \right) + \gd - \frac{\gd+1}{p} \right) p^{-(\gd/2+1)} \right| \leq \frac{2}{p},
\]
which can be proved by using elementary calculus techniques to bound the function
\[
y \mapsto \frac{2}{p} - \frac{2y}{p^{y/2 + 1/2}} + \frac{2y}{p^{y/2 + 3/2}} \quad\quad \left( y \in [1,\infty) \right).
\]
To prove \eqref{exactformula}, we reason as follows.  Writing a matrix of trace $r$ as
\[
A = 	\begin{pmatrix}
		a & b \\
		c & r-a
	\end{pmatrix},
\]
we are led to ask how many solutions $(a,b,c)$ modulo $p^n$ there are to the equation 
\begin{equation} \label{mainquadratic}
\left( a-\frac{r}{2} \right)^2 \equiv \frac{\gD -4bc}{4} \mod p^n,
\end{equation}
which leads us to the following two sublemmas, whose proofs are straightforward calculations.
\begin{sublemma}
For any odd prime power $p^n$, the number $N_y$ of solutions $x$ modulo $p^n$ to 
\[
x^2 \equiv y \mod p^n
\]
is given by
\[
N_y = 	\begin{cases}
			p^{n - \lceil n/2 \rceil} & \text{ if } y \equiv 0 \mod p^n \\
			p^m \left( 1 + \left( \frac{y/p^{2m}}{p} \right) \right) & \text{ if } y = p^{2m}\cdot y' \text{ with } p \nmid y' \text{ and } 2m < n \\
			0 & \text{ otherwise.}
		\end{cases}
\]
\end{sublemma}
\begin{sublemma}
For any odd prime power $p^n$, the number $P_{y}$ of pairs $(b,c)$ modulo $p^n$ satisfying
\[
4bc \equiv y \mod p^n
\]
is given by
\[
P_y = 	\begin{cases}
			(m+1) \cdot \varphi(p^n) & \text{ if } y = p^m \cdot y' \text{ with } p \nmid y' \text{ and } m < n \\
			n \cdot \varphi(p^n) + p^n & \text{ if } y \equiv 0 \mod p^n.
		\end{cases}
\]
\end{sublemma}
The number of solutions $(a,b,c)$ modulo $p^n$ to \eqref{mainquadratic} is simply
\[
\sum_{y \mod p^n} N_{\gD-y} \cdot P_y.
\]
Using the two sublemmas and some calculation, we arrive at \eqref{exactformula}, proving Lemma \ref{bonnlemma}.
\hfill $\Box$ \\

By Lemma \ref{bonnlemma}, we see that
\[
| \det{}^{-1}(\gd^{-1}(h))_r | \ll \varphi(m_2) \cdot m_2^2 \cdot \prod_{p \leq m_2} \left( 1 + \frac{1}{p} \right)^3 \leq m_2^3 \log^3 m_2.
\]
Thus we have
\[
\frac{m_E | G(m_E)_r |}{| G(m_E) |} \ll \log^3 m_2,
\]
which proves the bound for $C_{E,r}$ in the non-CM case. \\

For the Koblitz constant in the CM case, we see that
\[
C_{E,\text{prime}} \; \ll \; \frac{1}{\prod_{\ell \mid \gD_E}(1-1/\ell)} \; \ll \; \log(\gD_E).
\]
For $C_{E,\text{prime}}$ in the non-CM case, we similarly have
\[
C_{E,\text{prime}} \; \ll \; \frac{1}{\prod_{\ell \leq m_E}\left( 1 - 1/\ell \right) } \; \ll \; \log m_E,
\]
finishing the proof of Lemma \ref{constantbounds}.
\end{proof}
We now use the following result, which is a restatement of \cite[Theorem $3$]{jonestc}.
\begin{theorem}
Assume an affirmative answer to Question \ref{serresquestion}.  Then for any non-CM elliptic curve $E$ over $\mbq$ we have
\[
m_E \ll \left( \prod_{p \mid \gD_E} p \right)^5,
\]
with an absolute constant.
\end{theorem}
Note that, for $E \in \mc{C}$, we have
\[
\prod_{p \mid \gD_E} p \leq \gD_E \ll 4A^3+27B^2,
\]
and thus, we have
\[
\frac{1}{| \mc{C} |} \sum_{{\begin{substack} {E \in \mc{C} \\ E \text{ is not a Serre curve}} \end{substack}}} \left| C_E - C \right|^k \quad \ll_k \quad  \frac{ \log^{3k}(A^3+B^2) }{| \mc{C} |} \sum_{{\begin{substack} {E \in \mc{C} \\ E \text{ is not a Serre curve}} \end{substack}}} 1
 \]
We finally use the following result of \cite{jonessc}, which in our situation may be stated as follows:
\begin{theorem}
There is a $\gamma > 0$ so that
\[
\sum_{{\begin{substack} {E \in \mc{C} \\ E \text{ is not a Serre curve}} \end{substack}}} 1 \quad \ll \quad \frac{|\mc{C}| \log^\gamma ( \min\{A,B\} )}{\sqrt{\min\{A,B\}}},
\]
with an absolute implied constant.
\end{theorem}
This implies \eqref{overnonserrecurves}, finishing the proof of Theorem \ref{maintheoremcond}.

\end{document}